\documentclass[11pt]{amsart}

\usepackage{amsfonts}
\usepackage{verbatim, amsmath, amssymb}

\newtheorem{theorem}{Theorem}[section]
\newtheorem{question}[theorem]{Question}
\newtheorem{conjecture}[theorem]{Conjecture}

\theoremstyle{definition}
\newtheorem{example}[theorem]{Example}

\newcommand{\de}{\partial}
\newcommand{\db}{\overline{\partial}}

\newcommand{\ddbar}{i \partial \overline{\partial}}
\newcommand{\Ric}{\mathrm{Ric}}
\newcommand{\ov}[1]{\overline{#1}}

\newcommand{\ti}[1]{\tilde{#1}}
\newcommand{\vp}{\varphi}
\newcommand{\ve}{\varepsilon}

\renewcommand{\leq}{\leqslant}
\renewcommand{\geq}{\geqslant}

\begin{document}

\title{Ricci-flat metrics on Calabi-Yau manifolds}

    \author{Valentino Tosatti}

\address{Courant Institute of Mathematical Sciences, New York University, 251 Mercer St, New York, NY 10012}
\email{tosatti@cims.nyu.edu}

\begin{abstract}
We study the space of Ricci-flat K\"ahler metrics on a given Calabi-Yau manifold, pose a number of questions about their possible degenerations, and survey some recent results on these questions.
\end{abstract}	
\maketitle

\section{Introduction}
\subsection{K\"ahler metrics and Ricci curvature}
Let $X$ be a compact complex manifold, with $\dim_{\mathbb{C}}X=n$. Given a Hermitian metric $g$ on $X$, we define a real $2$-form $\omega$ by
$$\omega(V,W):=g(JV,W),\quad V,W\in TX,$$
where $J:TX\to TX$ is the complex structure. The $2$-form $\omega$ is of type $(1,1)$, in the sense that $\omega(JV,JW)=\omega(V,W)$. The metric $g$ is called K\"ahler if we have $d\omega=0$. Since $g$ can be obtained from $\omega$ as $g(V,W)=\omega(V,JW)$, we will often use $g$ and $\omega$ interchangeably, and the complex structure $J$ will be viewed as fixed. From now on, we will assume that our complex manifold $X$ admits a K\"ahler metric.

The de Rham cohomology group $H^2(X,\mathbb{R})$ contains the subspace $H^{1,1}(X,\mathbb{R})$ of de Rham classes of closed real $(1,1)$-forms, and the fundamental $\de\db$-Lemma of Kodaira shows that if $\alpha,\beta$ are two closed real $(1,1)$-forms with $[\alpha]=[\beta]$ in $H^2(X,\mathbb{R})$, then there is a smooth ``potential'' function $\vp:X\to\mathbb{R}$, unique up to an additive constant, such that $\alpha=\beta+\ddbar\vp$, where $d=\de+\db$ is the usual decomposition of the exterior derivative on a complex manifold. In local holomorphic coordinates $\{z_1,\dots,z_n\}$, we can write
$$\ddbar\vp=i\sum_{j,k=1}^n \frac{\de^2\vp}{\de z_j\de\ov{z}_k}dz_j\wedge d\ov{z}_k.$$

As a matter of notation, when we will write $[\alpha]\in H^{1,1}(X,\mathbb{R})$, we will understand that we have chosen a specific closed real $(1,1)$-form $\alpha$ representing this cohomology class. In local holomorphic coordinates we can write
$$\alpha=i\sum_{j,k=1}^n \alpha_{j\ov{k}}(z)dz_j\wedge d\ov{z}_k,$$
where $(\alpha_{j\ov{k}}(z))$ is a Hermitian matrix for all $z$ in the chart. We will write $\alpha\geq 0$ (resp. $>0$) whenever this matrix is semipositive definite (resp. positive definite) at every point (this notion is independent of the chosen local chart).

We will denote by
$$\mathcal{C}:=\{[\omega]\ |\ \omega\ \text{K\"ahler metric on }X\}\subset H^{1,1}(X,\mathbb{R}),$$
the subset of cohomology classes of K\"ahler metrics. It is easy to see that $\mathcal{C}$ is an open convex cone in $H^{1,1}(X,\mathbb{R})$, known as the K\"ahler cone.

Its closure $\ov{\mathcal{C}}$ (with respect to the Euclidean topology) is known as the nef cone. Given a closed real $(1,1)$-form $\alpha$, its cohomology class $[\alpha]\in H^{1,1}(X,\mathbb{R})$ belongs to $\ov{\mathcal{C}}$ if and only if for every $\ve>0$ there is $\vp_\ve\in C^\infty(X,\mathbb{R})$ such that $\alpha+\ddbar\vp_\ve+\ve\omega\geq 0,$ where $\omega$ is any fixed K\"ahler metric on $X$.

Given a K\"ahler metric $\omega$ on $X$, its Ricci curvature form $\Ric(\omega)$ (which encodes the same information as the Ricci curvature tensor of the Riemannian metric $g$) is given locally as
\begin{equation}\label{ric}
\Ric(\omega)=-\frac{i}{2\pi}\de\db\log\det(g).
\end{equation}
It is thus a closed real $(1,1)$-form, whose cohomology class is independent of $\omega$, since if $\ti{\omega}$ is another K\"ahler metric then
$$\Ric(\omega)-\Ric(\ti{\omega})=\frac{i}{2\pi}\de\db\log\frac{\det(\ti{g})}{\det(g)},$$
is an exact $(1,1)$-form (note that $\log\frac{\det(\ti{g})}{\det(g)}$ is a global smooth function on $X$). The cohomology class of $\Ric(\omega)$ thus only depends on the complex structure of $X$, and it is denoted by $c_1(X)\in H^{1,1}(X,\mathbb{R})$, the first Chern class of $X$. 

The set of all K\"ahler metrics on $X$ will be denoted by $\mathcal{H}$.
Given $[\alpha]\in\mathcal{C}$ a K\"ahler class, we denote by $\mathcal{H}_{[\alpha]}\subset\mathcal{H}$ the set of K\"ahler metrics $\omega$ on $X$ with $[\omega]=[\alpha]$.
We thus obtain a map
\begin{equation}\label{bij}
\mathrm{Ric}:\mathcal{H}_{[\alpha]}\to\{\psi\ |\ \psi\text{ closed real }(1,1)\text{-form, }[\psi]=c_1(X)\},\quad \omega\mapsto\Ric(\omega).
\end{equation}
The following fundamental result is the starting point of our story:
\begin{theorem}[Calabi-Yau Theorem \cite{Ca,Ya}]\label{cy}
If $(X,\alpha)$ is a compact K\"ahler manifold, then the map in \eqref{bij} is a bijection.
\end{theorem}
Injectivity was proved by Calabi in 1954 in \cite{Ca} (using a by now standard strong maximum principle argument), and surjectivity by Yau in 1976 in \cite{Ya}. To prove surjectivity, given a closed real $(1,1)$-form $\psi$ on $X$ cohomologous to $c_1(X)$, the $\de\db$-Lemma allows us to find a smooth function $F$ such that $\Ric(\alpha)=\psi+\frac{i}{2\pi}\de\db F$, and adding a constant to $F$ we may assume that $\int_Xe^F\alpha^n=\int_X\alpha^n$. Then Yau shows in \cite{Ya} that there exists a smooth function $\vp$ with $\alpha+\ddbar\vp>0$ a new K\"ahler metric on $X$ cohomologous to $\alpha$, solving the PDE
\begin{equation}\label{maa}
(\alpha+\ddbar\vp)^n=e^F\alpha^n.
\end{equation}
In local holomorphic coordinates, \eqref{maa} becomes the complex Monge-Amp\`ere equation
\begin{equation}\label{maa2}
\det\left(\alpha_{j\ov{k}}+\frac{\de^2\vp}{\de z_j\de\ov{z}_k}\right)=e^F\det(\alpha_{j\ov{k}}),
\end{equation}
and so recalling \eqref{ric}, we see that
$$\Ric(\alpha+\ddbar\vp)=\Ric(\alpha)-\frac{i}{2\pi}\de\db F=\psi,$$
which shows that the map in \eqref{bij} is indeed surjective.

\subsection{Calabi-Yau manifolds and Ricci-flat metrics}
Of special interest are the compact K\"ahler manifolds for which $c_1(X)=0$ in $H^{1,1}(X,\mathbb{R})$, so that the $(1,1)$-form which is identically zero belongs to the image of \eqref{bij}.
Such manifolds are called {\em Calabi-Yau}, and by Theorem \ref{cy} these are precisely the compact K\"ahler manifolds which admit a Ricci-flat K\"ahler metric (i.e. a K\"ahler metric $\omega$ with $\Ric(\omega)\equiv 0$).

If $X$ is Calabi-Yau, we have a map
\begin{equation}
\mathrm{CY}:\mathcal{C}\to \mathcal{H},
\end{equation}
which to a K\"ahler class $[\alpha]$ associates the unique Ricci-flat K\"ahler metric $\omega$ on $X$ with $[\omega]=[\alpha]$. This is an injective map, with image equal to the set of all Ricci-flat K\"ahler metrics on $X$. If we equip $\mathcal{C}\subset H^{1,1}(X,\mathbb{R})$ with the Euclidean topology, and $\mathcal{H}$ with the topology of smooth convergence, then the map ${\rm CY}$ is a homeomorphism with its image.

\begin{example}
If $X=\mathbb{C}^n/\Lambda$ is a complex torus, then $X$ is Calabi-Yau, and given any $[\alpha]\in\mathcal{C}$, $\mathrm{CY}([\alpha])$ is the unique translation-invariant K\"ahler metric cohomologous to $[\alpha]$ (where $X$, viewed as a Lie group, acts on itself by translations). These K\"ahler metrics are flat, and these are essentially the only examples where $\mathrm{CY}([\alpha])$ can be described explicitly in closed form (more precisely, the only other such examples are finite quotients of tori by holomorphic free actions).
\end{example}

Given $[\alpha]\in\mathcal{C}$, the space $\mathcal{H}_{[\alpha]}$, equipped with the smooth topology, is noncompact. However, there is a classical way to embed it into a compact space, using the theory of closed positive currents, see e.g. \cite{Demb} for the basics about these objects. In a nutshell, currents can be thought of as differential forms with distributional coefficients, and as such there in an exterior derivative operator acting on them. One can then define cohomology groups using currents, which turn out to be equal to those defined using smooth forms. There is a notion of positivity for real $(1,1)$-currents, which is parallel to semipositivity of real $(1,1)$-forms.

Let then $\mathcal{T}$ be the space of all closed real positive $(1,1)$-currents $T$ on $X$, and $\mathcal{T}_{[\alpha]}$ the space of such currents with $[T]=[\alpha]$. These spaces of currents are equipped with the weak topology. Then the natural maps
$$\mathcal{H}_{[\alpha]}\to \mathcal{T}_{[\alpha]}, \quad \mathcal{H}\to\mathcal{T}$$
are injective and continuous, the space $\mathcal{T}_{[\alpha]}$ is compact, and so is $\bigcup_{[\alpha]\in K}\mathcal{T}_{[\alpha]}$ where $K\subset H^{1,1}(X,\mathbb{R})$ is any compact set.
In particular, we obtain a continuous injection
\begin{equation}\label{inj}
\mathcal{C}\overset{\mathrm{CY}}{\to} \mathcal{H}\hookrightarrow\mathcal{T}.
\end{equation}

\subsection{Degenerations of Ricci-flat metrics}
We would like to understand the possible degenerations of the Ricci-flat K\"ahler metrics on $X$, as we approach the boundary of the K\"ahler cone $\mathcal{C}$, and thus also understand geometrically such boundary classes. For this, the first natural question is:

\begin{question}\label{q1} Does there exist a continuous map
\begin{equation}
\ov{\mathrm{CY}}:\ov{\mathcal{C}}\to \mathcal{T},
\end{equation}
which extends the continuous injection in \eqref{inj}?
\end{question}

If such a continuous map $\ov{\mathrm{CY}}$ exists, then for every $[\alpha]\in \ov{\mathcal{C}}$ we have $[\ov{\mathrm{CY}}([\alpha])]=[\alpha]$, since this is true when $[\alpha]\in\mathcal{C}$. In other words, $\ov{\mathrm{CY}}([\alpha])$ would be a closed positive real $(1,1)$-current in the cohomology class $[\alpha]\in\ov{\mathcal{C}}$ with the property that whenever a sequence of K\"ahler classes $\{[\alpha_i]\}_{i\geq 0}\subset\mathcal{C}$ converges to $[\alpha]$, the corresponding Ricci-flat K\"ahler metrics $\mathrm{CY}([\alpha_i])$ will converge weakly to it.

\begin{example}\label{torus}
If $X=\mathbb{C}^n/\Lambda$ is a complex torus, then the map $\ov{\mathrm{CY}}$ exists, and it associates to $[\alpha]\in\ov{\mathcal{C}}$ its unique translation-invariant representative, which is a semipositive definite real $(1,1)$-form with nontrivial kernel.
\end{example}

Another natural question is:
\begin{question}\label{q2} Given $[\alpha]\in\de\mathcal{C}$, given $\{[\alpha_i]\}_{i\geq 0}\subset\mathcal{C}$ with $[\alpha_i]\to[\alpha]$ and given $\eta\in\mathcal{T}_{[\alpha]}$ with
\begin{equation}\label{lim}
\mathrm{CY}([\alpha_i])\to \eta,
\end{equation}
in the weak topology.
\begin{itemize}
\item[(a)] What is the regularity of the current $\eta$?
\item[(b)] Does the convergence in \eqref{lim} happen in stronger topologies, perhaps away from a small singular set?
\end{itemize}
\end{question}
We can understand these questions in a more explicit way as follows. It is well-known that the $\de\db$-Lemma also holds for currents, and so $\eta$ can be written as $\eta=\alpha+\ddbar\vp$ for some distribution $\vp$. The fact that $\eta$ is a positive current implies that $\vp$ is the distribution associated to a quasi-plurisubharmonic function $\vp:X\to\mathbb{R}\cup\{-\infty\}$ (meaning in local holomorphic charts on $X$ it can be written as the sum of a plurisubharmonic function and a smooth function). The function $\vp$ is unique up to an additive constant, and is at least in $L^1(X)$ (in fact, in all $L^p(X)$ with $p<\infty$). Then part (a) asks about the regularity of the function $\vp$. For example, do we have $\vp\in L^\infty(X)$? Or perhaps $\vp\in C^\infty(X\backslash V)$  for some proper closed analytic subvariety $V\subset X$?

As for part (b), we can write $\mathrm{CY}([\alpha_i])=\alpha_i+\ddbar\vp_i$ and choose the additive constants properly so that we have $\vp_i\to\vp$ in $L^1(X)$ (again, also in all $L^p(X), p<\infty$). Part (b) then asks: if we know that $\vp$ lies in a better space (from part (a)), can we show that the convergence $\vp_i\to\vp$ also happens in this topology?

\begin{example}
Again when $X=\mathbb{C}^n/\Lambda$ is a complex torus, the limiting currents $\ov{\mathrm{CY}}([\alpha])$ are of course smooth on $X$, and the convergence in \eqref{lim} happens in the smooth topology.
\end{example}

\section{The non-collapsing case}
\subsection{The boundary of the K\"ahler cone}
Points in $\de\mathcal{C}=\ov{\mathcal{C}}\backslash\mathcal{C}$ are nef $(1,1)$-classes which are not K\"ahler. The following fundamental result is the natural extension of the Nakai-Moishezon Theorem in algebraic geometry to the more general setting of compact K\"ahler manifolds:

\begin{theorem}[Demailly-P\u{a}un \cite{DP}]\label{DP}
A nef $(1,1)$-class $[\alpha]\in \ov{\mathcal{C}}$ belongs to $\de\mathcal{C}$ if and only if there is a positive-dimensional irreducible closed analytic subvariety $V\subset X$ with
$$\int_V\alpha^{\dim V}=0.$$
\end{theorem}
This is complemented by the following transcendental extension of a theorem of Nakamaye \cite{Na} in algebraic geometry:
\begin{theorem}[Collins-T. \cite{CT}]\label{ct}
Given a nef $(1,1)$-class $[\alpha]\in \ov{\mathcal{C}}$, the set
$$\mathrm{Null}([\alpha]):=\bigcup_{\int_V\alpha^{\dim V}=0}V,$$
is a closed analytic subvariety of $X$. Furthermore, there is a closed positive current $T$ in $[\alpha]$ with $T\geq \ve\omega$ on $X$ ($\ve>0$), and $T$ is a smooth K\"ahler metric on $X\backslash \mathrm{Null}([\alpha])$.
\end{theorem}
Theorem \ref{DP} thus says that given $[\alpha]\in\ov{\mathcal{C}}$, we have $\mathrm{Null}([\alpha])=\emptyset$ if and only if $[\alpha]\in\mathcal{C}$. Theorem \ref{ct} says that away from $\mathrm{Null}([\alpha])$, the class $[\alpha]$ ``behaves like'' a K\"ahler class.

\subsection{Nef and big classes}

A nef $(1,1)$-class $[\alpha]\in\ov{\mathcal{C}}$ is called big if $\int_X\alpha^n>0$, or equivalently if $\mathrm{Null}([\alpha])\neq X$. The cone $\mathcal{B}$ of nef and big classes
clearly satisfies $\mathcal{C}\subset \mathcal{B}\subset\ov{\mathcal{C}}$.

Observe that if $[\alpha]\in\mathcal{B}$ and $\{[\alpha_i]\}_{i\geq 0}\subset\mathcal{C}$ with $[\alpha_i]\to[\alpha]$, then the Ricci-flat metrics $\omega_i:=\mathrm{CY}([\alpha_i])$ are volume non-collapsed in the sense that
$$\mathrm{Vol}(X,\omega_i)=\frac{1}{n!}\int_X\omega_i^n=\frac{1}{n!}\int_X\alpha_i^n\overset{i\to+\infty}{\longrightarrow}\frac{1}{n!}\int_X\alpha^n>0.$$

The following theorem gives a positive answer to Questions \ref{q1} and \ref{q2} for nef and big classes:

\begin{theorem}[T. \cite{To}, Boucksom-Eyssidieux-Guedj-Zeriahi \cite{BEGZ}, Collins-T. \cite{CT}]\label{deg}
There is a continuous injective map
$$\ov{\mathrm{CY}}:\mathcal{B}\to \mathcal{T},$$
which extends the continuous injection in \eqref{inj}. Furthermore, given $[\alpha]\in\mathcal{B}$ and $\{[\alpha_i]\}_{i\geq 0}\subset\mathcal{C}$ with $[\alpha_i]\to[\alpha]$ as $i\to\infty$, we have that $\ov{\mathrm{CY}}([\alpha])$ is a smooth Ricci-flat K\"ahler metric on $X\backslash \mathrm{Null}([\alpha])$ and
$$\mathrm{CY}([\alpha_i])\to\ov{\mathrm{CY}}([\alpha])$$
in $C^\infty_{\rm loc}(X\backslash \mathrm{Null}([\alpha]))$.
\end{theorem}
It is easy to see that $X\backslash \mathrm{Null}([\alpha])$ is the largest open set where this smooth convergence holds. Theorem \ref{deg} was first established by the author in \cite{To} when $X$ is projective and $[\alpha]$ belongs to the real N\'eron-Severi group, using also the work of Eyssidieux-Guedj-Zeriahi \cite{EGZ}. The existence of the continuous injective map $\ov{\mathrm{CY}}$ in the general setting of Theorem \ref{deg} follows from the later work of Boucksom-Eyssidieux-Guedj-Zeriahi \cite{BEGZ}, and the smooth convergence away from $\mathrm{Null}([\alpha]))$ is proved in \cite{CT} by Collins and the author.

\begin{example}
Consider a $2$-torus $\mathbb{C}^2/\Lambda$, and take its quotient $Y$ by the involution $\iota(z_1,z_2)=(-z_1,-z_2)$. Then $Y$ has $16$ singular points which can be resolved by blowing them up, and obtaining $\pi:X\to Y$ where $X$ is a $K3$ surface, known as the Kummer surface associated to the torus. Pulling back a K\"ahler class on $Y$, one obtains a nef and big class $[\alpha]$ on $X$ with $\mathrm{Null}([\alpha])=\mathrm{Exc}(\pi)$ consisting of $16$ disjoint rational curves. Then $\ov{\mathrm{CY}}([\alpha])$ is the pullback of a flat orbifold K\"ahler metric on $Y$. In this case, Theorem \ref{deg} was first proved by Kobayashi-Todorov \cite{KT} and LeBrun-Singer \cite{LS} via gluing methods.
\end{example}

\begin{example}
More generally, if $Y$ is a normal analytic variety with canonical singularities and torsion canonical divisor (``Calabi-Yau variety''), which admits a crepant resolution $\pi:X\to Y$ (i.e. $X$ is smooth and Calabi-Yau), then pulling back K\"ahler classes on $Y$ one obtains many further examples where Theorem \ref{deg} applies. In this case, $\ov{\mathrm{CY}}([\alpha])$ is the pullback of a ``singular Ricci-flat metric'' on $Y$ constructed by Eyssidieux-Guedj-Zeriahi in \cite{EGZ}.
\end{example}

\begin{question}
Given $X$ Calabi-Yau and $[\alpha]\in\mathcal{B}$ a nef and big class, write $\ov{\mathrm{CY}}([\alpha])=\alpha+\ddbar\vp$ with $\vp\in L^1(X)$. Do we have $\vp\in L^\infty(X)$? $\vp\in C^0(X)$? Or even $\vp\in C^\gamma(X)$ for some $\gamma>0$?
\end{question}

It is proved in \cite{BEGZ} that $\ov{\mathrm{CY}}([\alpha])$ has ``minimal singularities'' in the class $[\alpha]$, which means that if we write $\ov{\mathrm{CY}}([\alpha])=\alpha+\ddbar\vp$ and if $T=\alpha+\ddbar\psi$ is any other closed positive $(1,1)$-current in the class $[\alpha]$, then there is $C>0$ such that
\begin{equation}\label{min}
\psi\leq \vp+C,
\end{equation}
on $X$. The intuition here is that quasi-plurisubharmonic functions may equal $-\infty$ somewhere, and \eqref{min} says that $\vp$ goes to $-\infty$ at a slower rate than all other possible candidates (if at all).

Because of this, in order to prove that $\vp\in L^\infty(X)$ it would be sufficient to show that $[\alpha]$ contains at least one closed positive $(1,1)$-current with bounded local potentials. This would follow from the following much stronger conjecture:
\begin{conjecture}[Filip-T., Conj. 1.2 in \cite{FT}]\label{fi}
Given $X$ Calabi-Yau and $[\alpha]\in\mathcal{B}$ a nef and big $(1,1)$-class, there is $\pi:X\to Y$ a bimeromorphic holomorphic map onto a Calabi-Yau variety with $[\alpha]=\pi^*[\beta]$ for some K\"ahler class $[\beta]$ on $Y$. In particular, $[\alpha]$ contains a smooth semipositive form.
\end{conjecture}
This conjecture, which is a transcendental version of Kawamata's base-point-free theorem, is known when $n\leq 3$ by \cite{FT,Ho}, and for all $n$ when $X$ is projective \cite{DH}. Furthermore, if $[\alpha]=\pi^*[\beta]$ as above, then it is expected that $\vp$ should  be H\"older continuous, see e.g. \cite{CGZ,DZ,DV,GGZ} for works that establish that $\vp\in C^0(X)$ (or even log continuous) under various assumptions, and \cite{GKSS} for a recent result on H\"older continuity.

Another related conjecture is the following:
\begin{question}
In the setting of Theorem \ref{deg}, write $\mathrm{CY}([\alpha_i])=\alpha_i+\ddbar\vp_i,  \ov{\mathrm{CY}}([\alpha])=\alpha+\ddbar\vp,$ with $\vp_i\to\vp$ in $L^1(X)$. Do we have $\vp_i\to \vp$ in $L^\infty(X)$?
\end{question}
This would of course imply that $\vp\in C^0(X)$. The answer to this question is affirmative if $[\alpha]=\pi^*[\beta]$ where $\pi:X\to Y$ is bimeromorphic holomorphic and $Y$ is a projective Calabi-Yau variety and $[\beta]=c_1(L)$ where $L\to Y$ is an ample line bundle, thanks to \cite{EGZ} and \cite{CGZ}.

\section{The collapsing case}
\subsection{Nef classes which are not big}

Given the positive results in Theorem \ref{deg} for nef and big classes, it may come as a surprise that the answer to Question \ref{q1} is in fact negative in general:

\begin{theorem}[Filip-T. \cite{FT3}]\label{fil}
There exists a projective $K3$ surface $X$ with $[\alpha]\in \ov{\mathcal{C}}\backslash\mathcal{B}$ with two sequences $\{[\alpha_i]\}_{i\geq 0},\{[\beta_i]\}_{i\geq 0}\subset\mathcal{C}$ with $[\alpha_i]\to [\alpha],$ $[\beta_i]\to[\alpha]$ is $i\to\infty$, and $\eta_1,\eta_2\in\mathcal{T}_{[\alpha]}$ with $\eta_1\neq \eta_2$, such that
$$\mathrm{CY}([\alpha_i])\to \eta_1,\quad \mathrm{CY}([\beta_i])\to\eta_2,$$
in the weak topology of currents.
\end{theorem}

To construct $X$, take a pencil $\{C_t\}_{t\in\mathbb{P}^1}$ of cubic curves in $\mathbb{P}^2$ with general member smooth, and singular members with only nodal singularities. Choose $t\neq s$ in $\mathbb{P}^1$, and let $X$ be the double cover of the blowup of $\mathbb{P}^2$ at the $9$ base points of the pencil, ramified over the strict transforms of $C_t$ and $C_s$. This is a smooth projective $K3$ surface, with an elliptic fibration $\pi:X\to \mathbb{P}^1$. We then take $[\alpha]=\pi^*[\omega_{\rm FS}]$ and $[\alpha_i]=[\alpha]+\frac{1}{i}[\omega]$, where $\omega$ is a Ricci-flat K\"ahler metric on $X$. Then it is proved in \cite{To4,GTZ,HT} that $\mathrm{CY}([\alpha_i])\to \eta_1$ weakly and also smoothly on $X\backslash V$ (where $V$ are the singular fibers of $\pi$), where $\eta_1=c\pi^*\pi_*(\omega^2), c\in\mathbb{R}_{>0},$ and $\pi_*(\omega^2)$ equals the pushforward of the Calabi-Yau volume form $\omega^2$. Clearly, $\eta_1$ is smooth away from the singular fibers of $\pi$.

The elliptic fibration $\pi$ has at least $9$ sections, corresponding to the exceptional divisors of the blowups. Using one of these as the origin makes the smooth fibers into elliptic curves, and translating fiberwise to each of the other sections we obtain an action of $\mathbb{Z}^8$ on $X$ by automorphisms which preserve the elliptic fibration. Pick $T:X\to X$ one such automorphism with infinite order, and for $i\geq 1$ let
$$[\beta_i]:=c'\frac{(T^i)^*[\omega]}{i^2},$$
then if $c'>0$ is chosen appropriately, we have that $[\beta_i]\to [\alpha]$, see \cite{DF}. The weak limit
$$\mathrm{CY}([\beta_i])\to\eta_2,$$
exists, thanks to \cite[Cor. 3.2.20]{FT3}, and $\eta_2$ is a closed positive current in $[\alpha]$, which is necessarily of the form $\eta_2=\pi^*\gamma$ for some closed positive current $\gamma$ on $\mathbb{P}^1$, and $\eta_2$ is smooth away from the singular fibers of $\pi$. Now, in \cite[\S 4.4.12]{FT3} it is shown that the current $\gamma$ is independent of the choice of $s,t$ that we made at the beginning of the construction, while the current $\pi_*(\omega)^2$ depends nontrivially on the choice of $s,t$. Hence there exist values of $s,t$ for which $\eta_1=c\pi^*(\pi_*(\omega)^2)$ is different from $\eta_2=\pi^*\gamma$, as claimed.

As for question \ref{q2}, as we just mentioned, in this example we have that both $\eta_1$ and $\eta_2$ are smooth away from some proper closed analytic subvariety $V\subset X$, which settles part (a) of Question \ref{q2}. As for part (b), in this example we have that $\mathrm{CY}([\alpha_i])\to \eta_1$ in $C^\infty_{\rm loc}(X\backslash V)$, but interestingly it was observed by Cao \cite{Cao} that
the convergence $\mathrm{CY}([\beta_i])\to \eta_2$ does not in general happen in $C^0_{\rm loc}(X\backslash W)$ for any proper closed analytic subvariety $W\subset X$. This shows that Question \ref{q2} part (b) may fail.

\subsection{Collapsing setup}
Let $X$ be a Calabi-Yau manifold and $[\alpha]\in\ov{\mathcal{C}}$ a nef class which is not big, i.e.
$$\int_X\alpha^n=0.$$
We have seen in Theorem \ref{fil} above that in general we cannot hope to extend continuously the Calabi-Yau map $\mathcal{C}\to\mathcal{T}$ to the class $[\alpha]$, since approaching $[\alpha]$ in different ways with K\"ahler classes may produce different weak limits for the corresponding Ricci-flat K\"ahler metrics.

However, the situation is expected to improve if we only consider paths of the form
$$[\alpha]+t[\omega],$$
where $[\omega]\in\mathcal{C}$ and $t\in [0,1]$. Observe that $[\alpha]+t[\omega]\in\mathcal{C}$ for $t>0$. Then the Ricci-flat metrics
\begin{equation}\label{satan}
\omega_t=\mathrm{CY}([\alpha]+t[\omega]), \quad t\in (0,1],
\end{equation}
are volume collapsing as $t\downarrow 0$, since
$$\mathrm{Vol}(X,\omega_t)=\frac{1}{n!}\int_X\omega_t^n=\frac{1}{n!}\int_X(\alpha+t\omega)^n\overset{t\downarrow 0}{\longrightarrow} \frac{1}{n!}\int_X\alpha^n=0.$$
This setting was first studied explicitly by Gross-Wilson \cite{GW} when $X$ is a $K3$ surface with an elliptic fibration $\pi:X\to\mathbb{P}^1$ and $[\alpha]=\pi^*[\omega_{\rm FS}]$.

\begin{question}\label{q3}
Let $X$ be a Calabi-Yau manifold, $\omega$ a K\"ahler metric on $X$, and $[\alpha]\in\ov{\mathcal{C}}$ a nef class with $\int_X\alpha^n=0.$ For $t>0$ define Ricci-flat K\"ahler metrics $\omega_t$ on $X$ by \eqref{satan}.
\begin{itemize}
\item[(a)]
As $t\downarrow 0$, do $\omega_t$ converge in the weak topology to some closed positive current $\eta\in[\alpha]$?
\item[(b)]
Assuming (a), is $\eta$ independent of the choice of $[\omega]$?
\item[(c)] Assuming (a), do we have
$\omega_t\to\eta$ in $C^\infty_{\rm loc}(X\backslash V,g)$ for some proper closed analytic subset $V\subset X$?
\end{itemize}
\end{question}
About (a), since the set $K:=\{[\alpha]+t[\omega]\}_{t\in [0,1]}\subset H^{1,1}(X,\mathbb{R})$ is compact, the set $\bigcup_{[\beta]\in K}\mathcal{T}_{[\beta]}$ is compact in the weak topology, hence one can always extract some sequential weak limit $\eta$ of $\omega_{t_i}$, for some $t_i\downarrow 0$, but a priori the limit $\eta$ could depend on the sequence $t_i$. Part (c) is asking for uniform a priori $C^k$ estimates, for all $k\geq 0$, for the normalized solutions $\vp_t$ of the family of degenerating complex Monge-Amp\`ere equations
\begin{equation}\label{ma}
\omega_t^n=(\alpha+t\omega+\ddbar\vp_t)^n=c_t \omega^n,
\end{equation}
where $\omega$ can be assumed without loss to be Ricci-flat, and $c_t\in\mathbb{R}_{>0}$ is the constant given by
$$c_t=\frac{\int_X(\alpha+t\omega)^n}{\int_X\omega^n}\to 0.$$
Higher order estimates for $\omega_t$ (or equivalently for $\vp_t$) were conjectured by the author in 2010, see \cite{To2, To5}. In the case of elliptically fibered $K3$ surfaces with only $I_1$ singular fibers, the validity of Question \ref{q3} follows from the work of Gross-Wilson \cite{GW}, by constructing $\omega_t$ via a gluing method. Much more recently, gluing methods were employed to answer Question \ref{q3} for more elliptically fibered $K3$ surfaces \cite{CVZ} and for some Calabi-Yau $3$-folds with Lefschetz $K3$ fibrations \cite{Li}.

\subsection{Semiample classes}
Many nef classes which are not big arise from holomorphic fibrations, as follows: suppose $X$ is a Calabi-Yau manifold which admits a surjective holomorphic map $f:X^n\to Y^m$ with connected fibers onto a normal K\"ahler analytic space with $0\leq m<n$. If $\omega_Y$ is a K\"ahler metric on $Y$, then $\alpha=\pi^*\omega_Y$ is a smooth semipositive $(1,1)$-form on $X$ with $[\alpha]\in\ov{\mathcal{C}}$ and $\int_X\alpha^n=0$. In this case, there are proper closed analytic subvarieties $D\subset Y$ and $V:=f^{-1}(D)\subset X$ such that $Y\backslash D$ is smooth and $f|_{X\backslash V}:X\backslash V\to Y\backslash D$ is a proper holomorphic submersion. By the adjunction formula and Ehresmann's Lemma, its fibers $X_y=f^{-1}(y), y\in Y\backslash D,$ are all Calabi-Yau $(n-m)$-folds which are pairwise diffeomorphic, but not biholomorphic in general. We will refer to $V$ as the singular fibers.

Such classes $[\alpha]$ are called semiample, by analogy with algebraic geometry, since if $[\alpha]=c_1(L)$ and $L$ is semiample (in the usual sense that $L^k$ is base-point-free for some $k\geq 1$), then using sections of $L^m$ for $m$ sufficiently divisible we get a map $f:X\to Y$ as above. We can  rephrase Conjecture \ref{fi} above by saying that a nef and big class on a Calabi-Yau manifold should be semiample.

For nef classes which are not big, the following is a well-known conjecture in algebraic geometry (see e.g. \cite[Conj. 3.8]{To6}):

\begin{conjecture}
Let $X$ be a projective Calabi-Yau manifold, and $[\alpha]\in\ov{\mathcal{C}}$ a nef class with $\int_X\alpha^n=0$ and $[\alpha]\in H^2(X,\mathbb{Q})$. Then $[\alpha]$ is semiample.
\end{conjecture}

This remains open starting in dimension $3$. It is not hard to see that it holds when $X$ is a torus. Observe that a semiample class $[\alpha]$ need not belong to $H^2(X,\mathbb{Q})$ in general, even after rescaling it by any $\lambda\in\mathbb{R}_{>0}$. Also, in general nef classes which are not big are not necessarily semiample, see Section \ref{irrat} for an explicit example.

When $[\alpha]$ is semiample, so that it is the pullback of a K\"ahler class from the base of $f:X\to Y$ as above, Question \ref{q3} had received a lot of attention over the years, with progress made in \cite{Fi,ST,DPa,EGZ2,To4,GTZ,HT,TWY,TZ,HT2}. The final answer was obtained recently by Hein and the author \cite{HT3}:

\begin{theorem}[Hein-T. \cite{HT3}]\label{ht}
Suppose that $[\alpha]$ is semiample. Then Question \ref{q3} has a positive answer.
\end{theorem}

Parts (a) and (b) were proved early on by the author \cite{To4}, by showing that any weak limit $\eta$ has to be the pullback of a current from the base $Y$, which is of the form $\omega_Y+\ddbar\vp_0$ with $\vp_0$ smooth on $Y\backslash D$ and $\omega_Y+\ddbar\vp_0>0$ there, and which satisfies the complex Monge-Amp\`ere equation
\begin{equation}\label{mam}
(\omega_Y+\ddbar\vp_0)^{m}=c\pi_*(\omega_X^n),
\end{equation}
where $\omega_X$ is a fixed Ricci-flat K\"ahler metric on $X$ and $c>0$ is a constant (which does not depend on the choice of $\omega$). Since the solution of \eqref{mam} is known to be unique by \cite{EGZ,Ko}, this proves parts (a) and (b). On the other hand, part (c) turned out to be much harder. Here we give some ideas of the developments around this question, and of its resolution.
\subsubsection{No singular fibers}
If $Y$ is smooth and $f:X\to Y$ is a submersion everywhere, then global $C^\infty$ estimates for $\omega_t$ follow from the work of Fine \cite{Fi}, who constructed $\omega_t$ by deforming a semi-Ricci-flat form (see below) via the implicit function theorem. For this approach, the lack of singular fibers of $f$ is crucial. However, Zhang and the author have shown in \cite{TZ0,TZ2} that if $f$ has no singular fibers then it is necessarily a holomorphic fiber bundle (with $Y$ also Calabi-Yau), and if $X$ is projective it is in fact a finite quotient of a product. Because of this, the most interesting case is when the locus $V$ of singular fibers is nonempty.

\subsubsection{Quasi-isometry}
As a starting point, after normalizing $\vp_t$ with $\max_X\vp_t=0$ say, the uniform bound
\begin{equation}
\|\vp_t\|_{L^\infty(X)}\leq C,
\end{equation}
independent of $t\in (0,1]$ was proved by Demailly-Pali \cite{DPa} and Eyssidieux-Guedj-Zeriahi \cite{EGZ2} using pluripotential theory, and a new PDE proof was given recently by Guo-Phong-Tong \cite{GPT}.
This was then used as an ingredient by the author in \cite{To4}, and using also Yau's Schwarz Lemma computations \cite{Ya2}, to show that given any compact subset $K\subset X\backslash V$ there is $C>0$ such that for all $t\in (0,1]$ we have
\begin{equation}\label{dub}
C^{-1}(f^*\omega_Y+t\omega)\leq \omega_t\leq C(f^*\omega_Y+t\omega),
\end{equation}
on $K$. This means that on $K$, the metrics $\omega_t$ are uniformly quasi-isometric to the model shrinking metrics $f^*\omega_Y+t\omega$. In particular, the smooth fibers of $f$ are shrinking at the expected rate, and the ellipticity of the complex Monge-Amp\`ere equations \eqref{ma} is degenerating in the fiber directions. One does not expect any form of Schauder or Evans-Krylov theory to apply (uniformly in $t$) in such a degenerate limit, thus apparently ruling out the standard approach to higher order regularity.

\subsubsection{Uniform convergence}
The quasi-isometry statement in \eqref{dub} was later improved to a uniform convergence statement for $\omega_t$ by Weinkove, Yang and the author \cite{TWY}, as follows. Recall that the weak limit of $\omega_t$ is $\eta=f^*(\omega_Y+\ddbar\vp_0),$ which is smooth on $X\backslash V$. On the other hand, each smooth fiber $X_y=f^{-1}(y), y\in Y\backslash f(V)$, is a Calabi-Yau $(n-m)$-fold, so by the Calabi-Yau Theorem \ref{cy} there is $\rho_y\in C^\infty(X_y)$ with $\int_{X_y}\rho_y\omega^{n-m}=0$, with $\omega|_{X_y}+\ddbar\rho_y$ Ricci-flat K\"ahler on $X_y$. The functions $\rho_y$ vary smoothly in $y$, so taken together they define a smooth function $\rho\in C^\infty(X\backslash V)$, and we define
$$\omega_{\rm SRF}:=\omega+\ddbar\rho,$$
which is a closed real $(1,1)$-form on $X\backslash V$, called ``semi-Ricci-flat'' since $\omega_{\rm SRF}|_{X_y}$ is Ricci-flat, for all smooth fibers $X_y$. On $X\backslash V$ we can then write
$$\omega_t=\eta+t\omega_{\rm SRF}+\ddbar\psi_t,\quad \psi_t:=\vp_t-f^*\vp_0-t\rho.$$
The main result of \cite{TWY} is then that on any compact  subset $K\subset X\backslash V$  we have
\begin{equation}\label{exp1}
\omega_t=\eta+t\omega_{\rm SRF}+o(1),
\end{equation}
where the $o(1)$ term has $g_t$-norm that goes to zero on $K$ (which implies that its $g$-norm also goes to zero, thanks to \eqref{dub}). In particular, $\omega_t\to\eta$ in $C^0_{\rm loc}(X\backslash V,g)$.

\subsubsection{Higher order estimates in the case of torus fibers}
When the smooth fibers of $f$ are tori (or finite quotients of tori), Question \ref{q3} was solved affirmatively by Gross, Zhang and the author in \cite{GTZ}, assuming that $X$ is projective, and by Hein and the author \cite{HT} without this assumption. Here the rough idea is that the torus fibers can be ``unravelled'', passing to the universal cover of $f^{-1}(U)$ where $U\subset K$ is a sufficiently small open set. This universal cover turns out to be biholomorphic to $U\times\mathbb{C}^{n-m}$, via a fiber-preserving biholomorphism $p:U\times\mathbb{C}^{n-m}\to f^{-1}(U)$, and pulling back by the fiberwise stretching map $\lambda_t(y,z)=(y,t^{-\frac{1}{2}}z)$ on $U\times\mathbb{C}^{n-m}$ equation \eqref{dub} gives
\begin{equation}\label{dub2}
C^{-1}(p^*f^*\omega_Y+t\lambda_t^*p^*\omega)\leq \lambda_t^*p^*\omega_t\leq C(p^*f^*\omega_Y+t\lambda_t^*p^*\omega).
\end{equation}
To make the stretched PDE uniformly elliptic, it thus suffices to replace $\omega$ with a semipositive definite form $\omega_{\rm SF}\geq 0$ on $f^{-1}(U)$, which is positive definite in the fiber directions, such that
\begin{equation}\label{scale}
t\lambda_t^*p^*\omega_{\rm SF}=p^*\omega_{\rm SF}.
\end{equation}
This can indeed be done, as shown in \cite{GTZ,HT}, by taking
\begin{equation}\label{sf}
p^*\omega_{\rm SF}=\ddbar u,\quad u(y,z)=-\frac{1}{4}\sum_{i,j=1}^{n-m}\left(\mathrm{Im}Z(y)\right)^{-1}_{ij}(z_i-\ov{z}_i)(z_j-\ov{z}_j),
\end{equation}
where $Z:U\to\mathfrak{H}_{n-m}$ is a holomorphic period map (for the torus fibers) into the Siegel upper half space of symmetric $(n-m)\times (n-m)$ complex matrices with positive definite imaginary parts (that such a period map exists even without projectivity assumptions is shown in \cite{HT}). The form $\omega_{\rm SF}$ is ``semi-flat'' in the sense that it restricts to a flat K\"ahler metric on all the torus fibers over $U$. Remarkably, as shown in \cite[Thm 3.1]{TZ2}, the semi-flat form $\omega_{\rm SF}$ in \eqref{sf} agrees with the ``Betti form'' of this family of tori (when $X$ is projective), an object which has been much studied recently in Diophantine geometry, see e.g. \cite{ACZ}, where the scaling property \eqref{scale} also plays a crucial role.
In the context of elliptically fibered $K3$ surfaces, the semi-flat form was also used by Greene-Shapere-Vafa-Yau \cite{GSVY}, and its origin can be traced back to Satake's book \cite[Lemma IV.8.5]{Sat}.

Returning to our discussion, pulling back $\omega_t$ via $p \circ\lambda_t$ we get a Ricci-flat K\"ahler metric on $U\times\mathbb{C}^{n-m}$ which by \eqref{dub2} and \eqref{scale} is locally uniformly equivalent to Euclidean (i.e. the stretched PDE is now uniformly elliptic), and hence has uniform $C^k$ estimates for all $k\geq 0$ by Yau's estimates in \cite{Ya}. Stretching back, we obtain uniform $C^k$ estimates for our original $\omega_t$, thus proving part (c) of Question \ref{q3} in the case of torus fibers.

As observed in \cite[Appendix A]{DJZ} (see also \cite{JS}), the results of \cite{GTZ,HT} actually prove that
\begin{equation}\label{exp2}
\omega_t=\eta+t\omega_{\rm SRF}+\ddbar\underline{\psi_t}+o(1),
\end{equation}
where $o(1)$ goes to zero faster than $t^N$ for any $N\in\mathbb{N}$, together with all its derivatives, and where $\underline{\psi_t}$ is the pullback to $X\backslash V$ of the fiberwise average of $\psi_t$ with respect to the fiberwise Calabi-Yau volume form $(\omega_{\rm SRF}|_{X_y})^{n-m}$. Furthermore, the term $\ddbar\underline{\psi_t}$ goes to zero locally smoothly, albeit at an unknown rate.

\subsubsection{Higher order estimates in the isotrivial case}
Question \ref{q3} was solved affirmatively by Hein and the author \cite{HT2} when $f:X\to Y$ is isotrivial, which means that all smooth fibers $X_y$ are pairwise biholomorphic. By a classical theorem of Fischer-Grauert, $f:X\backslash V\to Y\backslash D$ is a holomorphic fiber bundle, so given any sufficiently small open set $U\subset Y\backslash D$ we can find a fiber-preserving biholomorphism between $f^{-1}(U)$ and $U\times F$, where $F=X_y$ is a fixed Calabi-Yau $(n-m)$-fold. In this trivialization our semi-Ricci-flat form $\omega_{\rm SRF}$ is equal to a Ricci-flat K\"ahler metric $\omega_F$ on the fiber $F$ trivially extended to the product $U\times F$.  The main result of \cite{HT2} is then that the same expansion in \eqref{exp2} holds, where again
$o(1)$ goes to zero faster than $t^N$ for any $N\in\mathbb{N}$, together with all its derivatives, and $\ddbar\underline{\psi_t}$ goes to zero locally smoothly. This of course implies that $\omega_t\to \eta$ in $C^\infty_{\rm loc}(X\backslash V,g)$. The main part of the proof of \eqref{exp2} is to show that for every $k\geq 1$, on $U\times F$ we have
\begin{equation}\label{bdd}
|\nabla^k \omega_t|^2_{g_t}\leq C_k,
\end{equation}
for all $t\in (0,1]$ (after shrinking $U$ slightly), where $\nabla$ is the covariant derivative of the product metric $\eta+\omega_F$. These bounds are proved by induction on $k$, and are obtained by a contradiction and blowup argument, relying also on Liouville theorem for Ricci-flat K\"ahler metrics proved in \cite{TZ,He,LLZ}. Roughly speaking, if the desired bound \eqref{bdd} fails to hold, so the supremum of $|\nabla^k \omega_t|^2_{g_t}$ blows up at some rate $\lambda_t^{2k}$ say ($\lambda_t\to+\infty$), we can stretch the base directions of $U\times F$ by $\lambda_t^{-1}$ (centered at a point where the supremum is roughly achieved), and rescale $\omega_t$ by $\lambda_t^{2}$, to obtain new Ricci-flat metrics $\ti{\omega}_t$ on say $B_{\lambda_t}\times F$ whose quantity in \eqref{bdd} is now bounded by $1$ (and equal to $1$ at some point), and which by \eqref{dub} are uniformly equivalent to $\omega_{\mathbb{C}^m}+\delta_t^2\omega_F$, where $\delta_t^2=t\lambda_t^2$. The argument then splits into $3$ cases according to whether we have a fast-forming blowup ($\delta_t\to+\infty$), regular-forming ($\delta_t\approx 1$) and slow-forming ($\delta_t\to 0$), up to passing to a sequence $t_i\to 0$ of course.

In the fast-forming case, zooming in by $\delta_t^{-1}$ near the blowup point, we obtain Ricci-flat metrics on $B_{\lambda_t}\times B_{\delta_t}$ which are uniformly quasi-isometric to Euclidean, hence by Yau's higher order estimates they will subsequentially converge to an entire solution, namely a Ricci-flat K\"ahler metric on $\mathbb{C}^n$ which is quasi-isometric to Euclidean, and is not constant. This however violates a well-known Liouville theorem \cite{RS}.

In the regular-forming case, the stretched and rescaled metrics $\ti{\omega}_t$ are already uniformly equivalent to a fixed metric, hence higher order estimates apply directly, and we can pass to the limit and obtain a Ricci-flat K\"ahler metric on $\mathbb{C}^m\times F$, quasi-isometric and cohomologous to a product metric, but not parallel with respect to it. Again, such a metric does not exist thanks to Hein's Liouville Theorem \cite{He, LLZ} (the earlier simpler version by Zhang and the author \cite{TZ} suffices here).

The slow-forming case requires much work. Here the stretched and rescaled metrics $\ti{\omega}_t$ are still collapsing in the fiber directions, and they subsequentially converge in $C^{k-1,\alpha}_{\rm loc}$ to the pullback of a K\"ahler metric on $\mathbb{C}^m$, which using arguments from \cite{To4} can be shown to be Ricci-flat, and quasi-isometric to Euclidean. By the Liouville theorem in $\mathbb{C}^m$, this metric must be constant. To derive a contradiction, we need to show that the limit is nonconstant, but this is not a priori obvious because of the collapsing of the fibers. Instead, we first prove that the stretched and rescaled metrics are locally uniformly bounded in $C^{k,\alpha}(g_t)$ (not just in $C^k(g_t)$, which is given by the blowup procedure), by linearizing and employing a new Schauder estimate on an infinite cylinder \cite[Thm. 3.14]{HT2} (proved using the method of Simon \cite{LSi}). This can then be used to show that in the norm $|\nabla^k \ti{\omega}_t|^2_{\ti{g}_t}$, which is equal to $1$ at the basepoint, the terms with all indices in the base directions are the dominant ones, and whenever a fiber direction appears the term is actually $o(1)$. Thus, passing to the limit, we see that the limiting metric is indeed nonconstant, which is a contradiction.

\subsubsection{The general case}
The method we just described, of using a blowup analysis and Liouville theorems, was refined in \cite{HT2} for general $f:X\to Y$, to prove that  $\omega_t\to\eta$ in $C^\gamma_{\rm loc}(X\backslash V,g)$, for any $\gamma<1$, and lastly in \cite{HT3} to upgrade this to $C^\infty_{\rm loc}(X\backslash V,g)$. The result proved there is in fact much more precise, and it provides an asymptotic expansion for $\omega_t$ generalizing \eqref{exp2} with uniform a priori estimates (which depend only on the constant $C$ in \eqref{dub}, and on the background data). To state these, we work on a sufficiently small ball $B$ inside a coordinate chart compactly contained in $Y\backslash D$, over which $f$ can be smoothly trivialized, so that there is a fiber-preserving diffeomorphism between $f^{-1}(B)$ and $B\times F$. The closed manifold $F$ is diffeomorphic to all the fibers $X_y, y\in B,$ and as such inherits a family of Ricci-flat Riemannian metrics $g_y, y\in B$, induced by the fiberwise Ricci-flat K\"ahler metric $\omega_{\rm SRF}|_{X_y}$ on $X_y$ via our diffeomorphism. It is important to realize that in general the metrics $g_y$ vary nontrivially in $y$, and in fact $g_y$ will not in general be parallel with respect to $g_{y'}$ for any $y'\neq y$ (unless we are in the cases described earlier, of isotrivial or torus fibers). The complex structure on $X$ pulled back via this diffeomorphism is a complex structure $J$ on $B\times F$, which is not a product structure in general, but for which the fibers $\{y\}\times F$ are $J$-holomorphic submanifolds. Recall that thanks to \eqref{dub}, on $U\times F$ the Ricci-flat metrics $g_t$ are uniformly quasi-isometric to the Riemannian product metrics $g_{y,t}:=g_{\mathbb{C}^m}+tg_y$ (which in general are not K\"ahler with respect to $J$). Taking inspiration from the isotrivial and torus-fibered cases, one might hope to prove that on $B\times F$ we have \eqref{exp2} where the $o(1)$ term has $\nabla^k$ with uniformly bounded $g_t$-norm for all $k\geq 1.$ But taking here $\nabla$ to be a product connection using one of the fiberwise Ricci-flat metrics will not work, since these norms depend very much on the choice of fiberwise metric, and do not remain uniformly equivalent to each other as $t\downarrow 0$. Rather, one is forced to define a connection $\mathbb{D}$ on $B\times F$ which along the fiber $\{y\}\times F$ (for any given $y\in B$) acts as the covariant derivative of $g_y$. This connection is torsion-free but it is not metric with respect to any fixed Riemannian metric, in general. Given thus $k\in\mathbb{N}$ and $0<\alpha<1$, taking $\mathbb{D}^k$, using the $\mathbb{D}$-parallel transport, and using $g_t$ to measure the tensor norms, we obtain a $t$-dependent family of ``shrinking H\"older norms'' on $B\times F$ that we will denote informally by $C^{k,\alpha}(g_t)$, and which are stronger than the standard $C^{k,\alpha}$ norm of any fixed Riemannian metric. The above results in the isotrivial and torus-fibered cases can then be restated by saying that \eqref{exp2} holds and the $o(1)$ term goes to zero in $C^{k,\alpha}(g_t)$ for all $k,\alpha$.

A major realization in \cite{HT3} was that in general this same statement is actually false, as long as $k\geq 2$, due to obstructions coming from the nontrivial variation of the complex structure of the fibers and from the non-flatness of the fiberwise Ricci-flat metrics $\omega_{\rm SRF}|_{X_y}$. With much work, we were able to identify concrete obstruction functions, which govern the higher-order terms in the expansion \eqref{exp2}, and were able to show that given any $k\geq 1$, on $B\times F$ we have
\begin{equation}\label{exp3}
\omega_t=\eta+t\omega_{\rm SRF}+\ddbar\underline{\psi_t}+\sum_{j=2}^k\gamma_{j,k}+\eta_{k},
\end{equation}
where $\|\eta_{k}\|_{C^{k,\alpha}(B\times F,g_t)}\leq C,$ $\|\gamma_{j,k}\|_{C^{k,\alpha}(B\times F, g)}\leq C,$
and $\ddbar\underline{\psi_t}$ goes to zero locally smoothly. Each term $\gamma_{j,k}$ for $2\leq j\leq k$, has $C^{j}(g_t)$-norm uniformly bounded, but its $C^{j,\alpha}(g_t)$-norm will blow up in general, although its (non-shrinking) $C^{k,\alpha}(g)$-norm remains bounded. Roughly speaking, one should think of $\gamma_{j,k}$ as being $\ddbar$ of $t^{\frac{j+2}{2}}$ times a fixed function on $B\times F$, although in reality the bounds we prove are slightly weaker than these (see \cite[Theorem 4.1]{HT3}). The decomposition in \eqref{exp3} is constructed in a highly nontrivial way from the above obstruction functions, and the proof of these {\em a priori} estimates for all the pieces is done by a very long contradiction and blowup argument in the same spirit as the above isotrivial case (but with many new difficulties). The first nontrivial piece $\gamma_{2,k}$ can be determined explicitly \cite[Theorem B]{HT3} and it equals
\begin{equation}
\gamma_{2,k}=t^2 \ddbar\Delta^{-1}_F \Delta^{-1}_F\left(\eta^{\mu\ov{\nu}}(\langle A_\mu, \ov{A_\nu}\rangle-\underline{\langle A_\mu, \ov{A_\nu}\rangle})\right),
\end{equation}
where $\eta$ is the limiting metric (viewed here as a metric on $B\subset Y\backslash D$), $\Delta_F$ is the fiberwise Laplacian on $\{y\}\times F$ with respect to the fiberwise Ricci-flat metric $g_y$, $\langle\cdot,\cdot\rangle$ is the fiberwise $g_y$-inner product, $A_\mu$ is the Kodaira-Spencer form in the direction $\de/\de z_{\mu}$ on $F$ (harmonic with respect to $g_y$), and underlining denotes the fiberwise average of a function with respect to $g_y$. Thus, $\gamma_{2,k}$ does not actually depend on $k$, and it is constructed purely in terms of the limiting metric $\eta$ and of the fiberwise Ricci-flat metrics. Observe that in the isotrivial case the Kodaira-Spencer forms vanish identically (hence $\gamma_{2,k}=0$), while in the torus-fibered case they do not vanish in general, but they are parallel and so $\langle A_\mu, \ov{A_\nu}\rangle=\underline{\langle A_\mu, \ov{A_\nu}\rangle}$ (hence again $\gamma_{2,k}=0$), which is consistent with \eqref{exp2}. These observation motivate the following:
\begin{question}\label{asym}
\begin{itemize}
\item[(a)] Can one prove a complete asymptotic expansion for $\omega_t$ where the terms do not depend on $k$?
\item[(b)] Can one find explicit formulas for the higher order terms in the expansion?
\item[(c)] Can one show that the term number $j$ in this expansion decays like $t^j$?
\end{itemize}
\end{question}
About (a), as mentioned above, the terms in the expansion \eqref{exp3} depend on the choice of $k$, but it is quite possible that one could show that each term converges to a limit when $k$ is sent to infinity, and that this would provide the desired $k$-independent expansion. About (b), in \cite{HT3} the term $\gamma_{2,k}$ is determined explicitly, and it is possible that the higher order terms could also be derived by arguing along those lines, but the combinatorial difficulties appear formidable. And about (c), the blowup arguments in \cite{HT3} so far fall just barely short of the decay of order $t^{\frac{j+2}{2}}$, which is well off the expected $t^j$ for $j>2$. It is tempting to imagine that microlocal analysis methods might help address Question \ref{asym}.

\subsubsection{K\"ahler-Ricci flow}
The results in \cite{HT3} that we just described were recently adapted by Hein, Lee and the author \cite{HLT} to a different setting. Here $X^n$ is a compact K\"ahler manifold with canonical bundle $K_X$ which is assumed to be semiample. We thus have a surjective holomorphic map $f:X\to Y$ with connected fibers onto a normal projective variety $Y$, with $c_1(K_X)=f^*c_1(L)$ for some ample line bundle $L\to Y$. The cases when $\dim Y=n$ and $\dim Y=0$ are by now well-understood, so we assume that $0<\dim Y<n$ (we refer to the author's recent survey \cite{To8} for more details). Given a K\"ahler metric $\omega_0$ on $Y$, let $\omega_t,t\geq 0,$ be the solution of the K\"ahler-Ricci flow
$$\frac{\de}{\de t}\omega_t=-\Ric(\omega_t)-\omega_t,$$
starting at $\omega_0$. It is known that this exists for all $t\geq 0$, and as $t\to+\infty$ the metrics $\omega_t$ exhibit an analogous collapsing behavior as the Ricci-flat metrics that we have just discussed. In particular, by proving a parabolic version of the above expansion, \cite{HLT} shows that $\omega_t$ collapse smoothly away from the singular fibers of $f$ to the pullback of a K\"ahler metric from $Y\backslash D$, and furthermore this collapsing happens with locally uniformly bounded Ricci curvature, thus resolving two conjectures of Song-Tian (see \cite{To8}).

\subsection{Non-semiample classes}\label{irrat}
As we have just explained in the previous section, when the nef class $[\alpha]$ is semiample, Question \ref{q3} has an affirmative answer. For general nef classes $[\alpha]$,
at the moment we do not have a positive answer to parts (a) and (b) in general, nor do we have a counterexample, and it would be very interesting to decide whether these hold or not. However, there are counterexamples to (c):

\begin{theorem}[Filip-T. \cite{FT}]\label{fil2}
There exist $K3$ surfaces $X$ with $[\alpha]\in\ov{\mathcal{C}}$ and $\int_X\alpha^n=0$, such that as $t\downarrow 0$ we have
$$\omega_t=\mathrm{CY}([\alpha]+t[\omega])\to \eta,$$
in the weak topology of currents, but not in $C^0_{\rm loc}(X\backslash V)$ for any proper closed analytic subvariety $V\subset X$. In fact, $\eta$ itself is not continuous on $X\backslash V$.
\end{theorem}

A very explicit example is obtained as follows: let $X=\{P=0\}\subset (\mathbb{P}^1)^3$ be a generic hypersurface of degree $(2,2,2)$. It is equipped with $3$ projections $X\to(\mathbb{P}^1)^2$, by forgetting one of the $\mathbb{P}^1$ factors, and these are ramified $2:1$ covers. We thus have holomorphic covering involutions $\sigma_i:X\to X, i=1,2,3,$ and their composition $T=\sigma_1\circ\sigma_2\circ\sigma_3$ is a holomorphic automorphism of $X$ with infinite order. The N\'eron-Severi group $NS(X)$ has rank $3$, spanned by the $3$ pullbacks of $\mathcal{O}_{\mathbb{P}^1}(1)$ under the $3$ projections $X\to\mathbb{P}^1$. With respect to this basis, the operator $T^*:NS(X)\to NS(X)$ is given by
\[\begin{pmatrix}
15 & 6 & 2\\
10 & 3 & 2\\
-6 & -2 & -1
\end{pmatrix}\]
whose largest eigenvalue is $\lambda=9+4\sqrt{5}$ with eigenvector the class $[\alpha]$ given by
\[[\alpha]=\left(\frac{18+8\sqrt{5}}{7+3\sqrt{5}},\frac{11+5\sqrt{5}}{7+3\sqrt{5}},-1\right).\]
The class $[\alpha]$ is nef with $\int_X\alpha^2=0$. A theorem of Gromov and Yomdin \cite{Gr} shows that the topological entropy of $T$ is $\log\lambda>0$, and a result of Cantat \cite{Can} then shows that $[\alpha]$ contains a unique closed positive current $\eta$. In particular, the Ricci-flat K\"ahler metrics $\omega_t=\mathrm{CY}([\alpha]+t[\omega])$ must converge weakly to $\eta$ as $t\downarrow 0$. If this convergence was uniform off a proper closed analytic subvariety $V\subset X$, then $\eta$ would be continuous on $X\backslash V$. However, a rigidity theorem in complex dynamics on projective surfaces, due to Cantat-Dupont \cite{CD}, shows that this happens only when $X$ is a Kummer surface, which is not the case here since Kummer surfaces have N\'eron-Severi rank at least $16$. Filip and the author gave a new proof of this result \cite{FT2} for $K3$ surfaces using Ricci-flat metrics, which also applied to non-projective $K3$s. Applying this result, one can also construct counterexamples to Question \ref{q3} (c) with $X$ one of the non-projective $K3$ surfaces constructed by McMullen \cite{McM} which have an automorphism with positive entropy whose dynamics admits a Siegel disc. See also the author's survey \cite{To7} on this topic.

\subsection{Weaker regularity}
Let us return to the general Question \ref{q2}. As we have seen in Theorem \ref{deg}, if the nef class $[\alpha]$ is also big, then the limiting current $\eta$ is smooth off a proper closed analytic subvariety, and the convergence $\mathrm{CY}([\alpha_i])\to \eta$ is also smooth there. And if the nef class $[\alpha]$ is not big, we have seen in Theorem \ref{fil2} that in general the weak limit $\eta$ of the Ricci-flat metrics $\mathrm{CY}([\alpha]+t[\omega])$ is not smooth (or even continuous) on a Zariski open set. The natural question is then whether some weaker regularity for $\eta$ would hold in general. In all of the above examples, the limiting current $\eta$ was of the form $\eta=\alpha+\ddbar\vp$ with $\vp\in L^\infty(X)$ (and often in $C^0(X)$ or even $C^\gamma(X)$ for some $\gamma<1$). Motivated by this, in \cite[Conjecture 3.7]{To3} the author posed the following:
\begin{conjecture}[T. \cite{To3}]\label{to}
Let $X$ be a Calabi-Yau manifold, and $[\alpha]\in\ov{\mathcal{C}}$ be a nef class. Then there is $\vp\in L^\infty(X)$ such that $T=\alpha+\ddbar\vp\geq 0$ is a closed positive current in $[\alpha]$.
\end{conjecture}
When $[\alpha]$ is not big, even the following weaker statement is open: $[\alpha]$ contains a closed positive current $T=\alpha+\ddbar\vp\geq 0$ with $\vp$ quasi-plurisubharmonic function with vanishing Lelong numbers $\nu(\vp,x)$, for all $x\in X$. Observe that when $[\alpha]$ is big and nef, Conjecture \ref{fi} implies Conjecture \ref{to} since it would even give us a smooth semipositive form in the class $[\alpha]$.

It is remarkable that for nef classes which are not big, Conjecture \ref{to} is open even on $K3$ surfaces. Nevertheless, Filip and the author obtained the following partial result in \cite{FT3}: if $X$ is a projective $K3$ surface with no $(-2)$-curves and with N\'eron-Severi group of rank at least $3$, then Conjecture \ref{to} holds for classes $[\alpha]\in\ov{\mathcal{C}}$ which belong to the real N\'eron-Severi group. In fact, the functions $\vp$ that we construct are not only bounded but continuous on $X$. These results exploit the dynamics of the $\mathrm{Aut}(X)$-action on $X$, and the Kawamata-Morrison cone conjecture, which is known for projective $K3$ surfaces \cite{Ste}.

To connect Conjecture \ref{to} to Question \ref{q2}, we pose another conjecture (recall the notion of minimal singularities in \eqref{min}):
\begin{conjecture}\label{to2}
Let $X$ be a Calabi-Yau manifold, $[\alpha]\in\ov{\mathcal{C}}$ be a nef class, and $\{[\alpha_i]\}_{i\geq 0}\subset\mathcal{C}$ with $[\alpha_i]\to[\alpha]$. Suppose that we have $\eta\in\mathcal{T}_{[\alpha]}$ with
\begin{equation}
\mathrm{CY}([\alpha_i])\to \eta,
\end{equation}
in the weak topology. Then $\eta$ has minimal singularities in the class $[\alpha]$.
\end{conjecture}

This conjecture is known when $[\alpha_i]=[\alpha]+\ve_i[\omega]$ with $\ve_i>0, \ve_i\to 0$, thanks to Fu-Guo-Song \cite{FGS} (see also \cite{GPTW} for a new proof), but is open for more general sequences. Whenever both Conjectures \ref{to} and \ref{to2} are known hold, we conclude that $\eta$ has bounded potentials, and furthermore that if we write $\mathrm{CY}([\alpha_i])=\alpha_i+\ddbar\vp_i$ (with $\max_X\vp_i=0$ say), then we have the uniform $L^\infty$ estimate
\begin{equation}\label{inf}
\|\vp_i\|_{L^\infty(X)}\leq C.
\end{equation}
Conversely, if one can prove \eqref{inf} then both Conjectures \ref{to} and \ref{to2} follow.

Lastly, we would like to pose the following conjecture about irrational nef but not big classes on $K3$ surfaces:
\begin{conjecture}\label{rig}
Let $X$ be a $K3$ surface, $[\alpha]\in\de\mathcal{C}$ with $\int_X\alpha^2=0$ and $\mathbb{R}[\alpha]\cap H^2(X,\mathbb{Q})=\{0\}$. Then $[\alpha]$ contains a unique closed positive $(1,1)$-current.
\end{conjecture}
Classes which contain only one closed positive current are called rigid. Sibony-Soldatenkov-Verbitsky \cite{SSV} recently introduced a method to prove such rigidity statements in many cases, see also the exposition by Filip and the author in \cite[Theorem 4.3.1]{FT3}. In particular, we know that Conjecture \ref{rig} holds when $X$ is projective with no $(-2)$-curves and with N\'eron-Severi group of rank at least $3$, and $[\alpha]$ is in the real N\'eron-Severi group. However, as far as the author can tell, this method has not yet settled Conjecture \ref{rig} in general, even assuming that $X$ is projective.

\section{Gromov-Hausdorff limits}

In this last section, we discuss limits of Ricci-flat metrics in a different topology: instead of considering weak limits of $\mathrm{CY}([\alpha_i])$ as currents, we regard $(X,\mathrm{CY}([\alpha_i]))$ as a compact metric space, and consider its limits in the Gromov-Hausdorff topology.

\subsection{Diameter}
First, we consider the behavior of the diameter of the Ricci-flat metrics. For this, a basic result is the following:
\begin{theorem}[T. \cite{To}, Zhang \cite{Zh}]\label{tz}
Let $X$ be a Calabi-Yau manifold and $K\subset H^{1,1}(X,\mathbb{R})$ be a compact subset. Then there is $C>0$ such that for all $[\omega]\in\mathcal{C}\cap K$ we have
$$\mathrm{diam}(X,\mathrm{CY}([\omega]))\leq C.$$
\end{theorem}
In particular, this implies a uniform diameter upper bound for Ricci-flat metrics of the form $\mathrm{CY}([\alpha_i])$ where $[\alpha_i]\in\mathcal{C}$ and $[\alpha_i]\to[\alpha]\in\de\mathcal{C}$. Interestingly, there are examples where the limit class $[\alpha]$ is not the zero class, and yet $\mathrm{diam}(X,\mathrm{CY}([\alpha_i]))\to 0$. Indeed, this happens for the Ricci-flat metrics in Theorem \ref{fil2}.

Even worse, in the examples in Theorem \ref{fil} we have
\begin{equation}\label{exag}
\mathrm{diam}(X,\mathrm{CY}([\alpha_i]))\geq c>0, \quad \mathrm{diam}(X,\mathrm{CY}([\beta_i]))\to 0,
\end{equation}
even though both sequences $[\alpha_i]$ and $[\beta_i]$ limit to the same class $[\alpha]$. This shows that the map
$$\mathcal{C}\to\mathbb{R}, \quad [\alpha]\mapsto \mathrm{diam}(X,\mathrm{CY}([\alpha])),$$
cannot be extended continuously to a map $\ov{\mathcal{C}}\to \mathbb{R}$.

Motivated by these examples, in \cite[Conjecture 3.25]{To3}, the author posed the following conjecture (indeed, the special case of segments in cohomology):
\begin{conjecture}\label{zer}
Let $X$ be a $K3$ surface, $[\alpha]\in\de\mathcal{C}$ with $\int_X\alpha^2=0$ and $\mathbb{R}[\alpha]\cap H^2(X,\mathbb{Q})=\{0\}$. Then for every sequence $[\alpha_i]\in\mathcal{C}$ with $[\alpha_i]\to[\alpha]$ we have
\begin{equation}\label{tozero}
\mathrm{diam}(X,\mathrm{CY}([\alpha_i]))\to 0.
\end{equation}
\end{conjecture}
Of course, \eqref{tozero} means that the Gromov-Hausdorff limit of $(X,\mathrm{CY}([\alpha_i]))$ is a point. 

\subsection{Gromov-Hausdorff limits}
Let $X$ be a Calabi-Yau manifold, with a sequence $[\alpha_i]\in\mathcal{C}$ with $[\alpha_i]\to[\alpha]\in\de\mathcal{C}.$
Thanks to the general diameter upper bound in Theorem \ref{tz}, Gromov's precompactness theorem shows that up to passing to a subsequence, we have that $(X,\mathrm{CY}([\alpha_i]))$ converges in the Gromov-Hausdorff topology to some compact metric space $(Z,d)$, which in general will  depend on the subsequence (recall \eqref{exag}). The main question is then:
\begin{question}\label{q4}
Can the Gromov-Hausdorff limit $(Z,d)$ be identified explicitly? When is it independent of the subsequence? Is it regular outside a singular set of high Hausdorff codimension?
\end{question}

\begin{example}
When $X$ is a torus, as explained in Example \ref{torus}, we have the smooth semipositive form $\ov{\mathrm{CY}}([\alpha])$, whose kernel defines a holomorphic foliation on $X$. The Gromov-Hausdorff limit of $(X,\mathrm{CY}([\alpha_i]))$ is then a flat real torus $T^{2n-k}$ where $k$ is the real dimension of the closure of any leaf of the foliation.
\end{example}

If $[\alpha]$ is nef and big, the answer is affirmative:
\begin{theorem}[Rong-Zhang \cite{RZ, CT}]
If $[\alpha]$ is nef and big, then $(X,\mathrm{CY}([\alpha_i]))$ converges in Gromov-Hausdorff to the metric completion of $(X\backslash\mathrm{Null}([\alpha]),\ov{\mathrm{CY}}([\alpha]))$.
\end{theorem}
If $[\alpha]=c_1(L)$ where $L$ is semiample, so there is a bimeromorphic holomorphic map $f:X\to Y$ and $L=f^*A$ with $A$ an ample line bundle over $Y$, then Song \cite{So} (improving results of Donaldson-Sun \cite{DS}) shows that $Z$ is homeomorphic to $Y$, and that $Z$ is regular outside a set of real Hausdorff codimension $4$. See also \cite{Sz2} for a new proof. In \cite[Conjecture 3.8]{To3} the author conjectured that even if $[\alpha]$ is not semiample, the Gromov-Hausdorff limit $(Z,d)$ should exist, be homeomorphic to a normal compact complex analytic space, and again have singularities in codimension $4$. The methods of \cite{So,Sz2} would likely give this if one could settle Conjecture \ref{fi}.

When $[\alpha]$ is nef but not big, the only known general results are when $[\alpha]$ is semiample, so we have a surjective holomorphic map $f:X\to Y$ with $\dim Y<n$ and $[\alpha]=f^*[\omega_Y]$, and we look at
$$\omega_t=\mathrm{CY}([\alpha]+t[\omega]), \quad t\in (0,1].$$
Recall that thanks to Theorem \ref{ht} we know that $\omega_t$ converges smoothly away from the singular fibers of $f$ to $\eta=f^*\omega_0$, where $\omega_0$ is a K\"ahler metric on $Y\backslash D.$
In this case, the following results were conjectured by the author in \cite{To2,To5}:
\begin{theorem}
In this setting, we have:
\begin{itemize}
\item[(a)] As $t\downarrow 0$, $(X,\omega_t)$ converges in Gromov-Hausdorff to $(Z,d)$, the metric completion of $(Y\backslash D, \omega_0)$.
\item[(b)] The metric space $(Z,d)$ is regular outside  a set of real Hausdorff codimension $2$.
\item[(c)] $Z$ is homeomorphic to $Y$.
\end{itemize}
\end{theorem}
This theorem was first proved when $\dim Y=1$ by Gross, Zhang and the author \cite{GTZ2}, when $X$ is hyperk\"ahler by Zhang and the author \cite{TZ2} (assuming $X$ projective, but this assumption is not needed because in this case $f$ is locally projective, see \cite[p.177]{LT2}), and when $Y$ is smooth and the divisorial part of $D$ has simple normal crossings by Gross, Zhang and the author \cite{GTZ3} (again the projectivity assumption there is unnecessary, see \cite[p.170]{LT}), see also \cite{Li2,Li3} for Lefschetz $K3$ fibrations and \cite{OO} for related work. Part (a) was later proved in general by Song-Tian-Zhang \cite{STZ}, who also proved part (c) if $Y$ is smooth (they assume $X$ projective, but it is straightforward to extend their results to the general case). Finally, parts (b) and (c) were recently settled by Sz\'ekelyhidi \cite{Sz} (he also assumes projectivity, to use the results of \cite{STZ}). 

To close, when $X$ is a $K3$ surface, a general result of Sun-Zhang \cite{SZ} gives a list of all possible Gromov-Hausdorff limits of all unit-diameter hyperk\"ahler metrics on $X$ (with respect to any complex structure), see also \cite{OY}. It is however likely that some of these limits (those which have odd real Hausdorff dimension) cannot actually arise as Gromov-Hausdorff limits of $(X,\mathrm{CY}([\alpha_i]))$ as above.

\section*{Acknowledgments}
The author is very grateful to all his co-authors for countless discussions on these topics over the years, and for the ideas they shared in our collaborations. He also thanks G\'abor Sz\'ekelyhidi and Junsheng Zhang for feedback on an earlier version of this paper. The author was partially supported by NSF grant DMS-2404599.

\end{document}